\documentclass[11pt]{amsart} 
\usepackage{graphicx}
\usepackage{amssymb, amsmath}
\usepackage{amsfonts}
\usepackage{amssymb}
\usepackage{float}

\newtheorem{theorem}{Theorem}

\newtheorem{definition}[theorem]{Definition}

\newtheorem{remark}[theorem]{Remark}

\newcommand{\divv}{\text{\rm div}}

\newcommand{\cL}{\mathcal L}
\def\cM{\stackrel{\circ}{M}}
\newcommand{\C}{\mathbb C}

\def\p{\partial}
\def\P{\partial}

\newcommand{\R}{\mathbb R}

\newcommand{\dist}{\text{\rm dist}}
\def\be{\begin{equation}}
\def\ee{\end{equation}}
\def\la{\lambda}
\def\TT{{\tau}}

\def\pr{{P_r}}

\def\R{\mathbb R}
\def\mA{\mathcal A}
\def\sR{\sqrt{R}}
\def\p{\partial}
\def\tilde{\widetilde}

\numberwithin{equation}{section}
\numberwithin{theorem}{section}
\numberwithin{figure}{section}

\begin{document}
\bibliographystyle{siam}

\title[Rayleigh-B\'enard Convection]{Rayleigh-B\'enard Convection: Dynamics and Structure in the Physical Space}

\author[T. Ma]{Tian Ma}
\address[TM]{Department of Mathematics, Sichuan University,
Chengdu, P. R. China and Department of Mathematics,
, Bloomington, IN 47405}

\author[S. Wang]{Shouhong Wang}
\address[SW]{Department of Mathematics,
Indiana University, Bloomington, IN 47405}
\email{showang@indiana.edu}

\thanks{The work was supported in part by the
Office of Naval Research,  by the National Science Foundation, 
and by the National Science Foundation
of China.}


\begin{abstract} The main objective of this article is part of a research program to link the dynamics of fluid flows with the structure and its transitions in the physical spaces. 
As a prototype of problem and to demonstrate the main ideas, we study the two-dimensional Rayleigh-B\'enard convection.  The analysis is based on two recently developed nonlinear theories: geometric theory for incompressible flows \cite{amsbook} and the bifurcation and stability theory for nonlinear dynamical systems (both finite and infinite dimensional) \cite{b-book}.
We have shown in \cite{benard} 
that the Rayleigh-B\'enard  problem bifurcates from the basic state to an 
attractor $\mA_R$ when the
Rayleigh number $R$ crosses the first critical Rayleigh number $R_c$ 
for all physically  sound boundary conditions, regardless of the 
multiplicity of the eigenvalue $R_c$ for the linear problem.
In this article, in addition to a classification of 
the bifurcated attractor $\mA_R$, the structure and its transitions of the solutions 
in the physical space is classified, leading to the existence and stability of two different flows structures: pure rolls and rolls separated by a cross the channel flow.  
It appears that the structure with rolls separated by a cross channel flow has not been carefully examined although it has been  observed in other physical contexts such as the Branstator-Kushnir waves in the atmospheric dynamics
\cite{branstator,kushnir}.
\end{abstract}
\maketitle

\section{Introduction}
\label{sc1}
The Rayleigh-B\'enard convection problem  was originated 
in the famous experiments conducted by H.\ B\'enard in 1900. 
B\'enard investigated a fluid, with a free surface, heated from 
below in a dish, and noticed a rather regular cellular pattern of hexagonal 
convection cells. Based on the pioneering studies by 
Lord Rayleigh \cite{rayleigh},  the convection
would occur only when the non-dimensional parameter, called the
Rayleigh number,
\begin{equation}
\label{eq1.1}
R = \frac{g\alpha\beta}{\kappa\nu}\, h^4
\end{equation}
exceeds a certain critical value, where $g$ is the acceleration due
to gravity, $\alpha$ the coefficient of thermal expansion of the
fluid, $\beta=|dT/dz|=(\bar T_0-\bar T_1)/h$ the vertical temperature 
gradient with $\bar T_0$ the temperature on the lower surface and $\bar T_1$ 
on the upper surface, $h$ the
depth of the layer of the fluid, $\kappa$ the thermal diffusivity and
$\nu$ the kinematic viscosity.
There have been intensive studies for this problem; see among others 
Chandrasekhar \cite{chandrasekhar} and Drazin and Reid \cite{dr} for 
linear theories, and   Kirchg{\"a}ssner \cite{kirch}, 
Rabinowitz \cite{rabinowitz}, and 
Yudovich \cite{yudovich67a,yudovich67b}, and the references therein for 
nonlinear theories.

Recently, the authors have developed a bifurcation theory \cite{b-book}
for nonlinear partial differential equations, which has been used to 
develop a nonlinear analysis for the Rayleigh-B\'enard convections 
\cite{benard}. 
This bifurcation theory is centered at a new notion of bifurcation, called 
attractor bifurcation for nonlinear 
evolution equations. The main ingredients of the theory include 
a) the attractor bifurcation theory, b) steady state bifurcation for a 
class of nonlinear problems with even order non-degenerate 
nonlinearities, regardless of the multiplicity of the 
eigenvalues, and c) new strategies for 
the Lyapunov-Schmidt reduction and the center manifold 
reduction procedures. 

In particular, based on this bifurcation theory, 
for the Rayleigh-B\'enard convection 
problem, we have shown \cite{benard, b-book} that 
the problem bifurcates from the trivial solution to an 
attractor $\mA_R$ when the
Rayleigh number $R$ crosses the first critical Rayleigh number $R_c$ 
for all physically  sound boundary conditions, regardless of the 
multiplicity of the eigenvalue $R_c$ for the linear problem. 

The main objectives of this article  are 1) to classify 
the solutions in the bifurcated attractor $\mA_R$, and 2) to study  
the structure and its transition of the solutions of the B\'enard problem 
in the physical space. The first objective is 
an important part of the above mentioned new bifurcation and stability theory. 
The second objective is part of a program 
recently initiated by the authors 
to develop a geometric theory 
of two-dimensional  
incompressible fluid flows in the physical spaces; see \cite{amsbook}.
This program of study  
consists of research in two directions: 
1) the study of the structure and its transitions/evolutions 
of divergence-free vector fields (kinematics), and 2) the study of 
the structure and its transitions 
of velocity fields for 2-D incompressible fluid flows governed 
by the Navier-Stokes equations or the Euler equations.
The study in this article is in the second direction, 
linking kinematics to dynamics.

To demonstrate ideas, in this article, we only consider 
the two-dimensional B\'enard convection problem.  The three-dimensional 
case is technically more involved, and shall be reported elsewhere. 
From the physical point of view, two-dimensional Boussinesq equations 
can be considered as idealized models for many physical phenomena, including 
1) the Walker circulation over the tropics \cite{salby}, 
which has the same topological structure as the cells given in 
Figure~\ref{fg4.3} in Theorem~\ref{th4.2}, 
 and 
2) the Branstator-Kushnir waves in the atmospheric dynamics 
\cite{branstator,kushnir}, which have similar topological 
structure as given in Figure~\ref{fg4.2} in Theorem~\ref{th4.1}. 

We end this introduction with a few remarks. 
First, the main idea of the study is to explicitly reduce the bifurcation problem 
to the center manifold, together with an $S^1$ attractor bifurcation theorem
and structural stability theorems for 2D incompressible flows. 
The types of solutions in this $S^1$ attractor depend on the boundary conditions. With the periodic boundary condition (\ref{eq2.7}) in the $x_1$ direction in this article, the bifurcated attractor consists of only steady states. When the boundary conditions for the velocity field  are free slip boundary conditions  and $\Omega=(0, L)^2\times (0, 1)$ with $0 < L^2 < (2-\sqrt[3]{2})/(\sqrt[3]{2}-1)$, using the same method proved in this article, we can prove that the bifurcated attractor is still an $S^1$, consisting of exactly eight singular steady states (with four saddles and four minimal attractors) and eight heteroclinic orbits connecting these steady states. The bifurcated attractor and its detailed classification provide a  global dynamic transitions in both the physical and phase spaces.
 
Second, the method and ideas presented in this article are crucial to obtain these results, which can not be obtained using only the classical bifurcation theories. For the case studied in this article, the classical bifurcation theory with symmetry arguments implies that  the bifurcation attractor in the main theorems, Theorems 4.1 and 4.2, {\it contain} a circle of steady states. We need, however,  the new bifurcated theory to prove in particular that the bifurcated attractors are {\it exactly} an $S^1$. Furthermore, for general boundary conditions such as the free-slip boundary conditions mentioned above, no symmetry can be used, and the classical amplitude equation methods fails to derive the dynamics.

Third, the newly developed geometric theory for incompressible flows is crucial for  the structure and its stability of the solutions in the physical spaces obtained in the main theorems. 

Fourth, it appears that the structure with rolls separated by a cross channel flow has not been carefully examined although it has been  observed in other physical contexts such as the Branstator-Kushnir waves in the atmospheric dynamics 
\cite{branstator,kushnir}.

Finally, as mentioned before, this article is part of a research program initiated in the mid 90s to make connections between the dynamics and the structure in the physical spaces. 

This article is organized as follows. First the functional setting and the attractor bifurcation theorem for the B\'enard convection 
problem obtained in \cite{benard} are introduced in Section 2.
Section 3 recapitulates 1) 
the approximation of the center manifold function, 
2) $S^1$ attractor bifurcation theorem, and  3) 
structural stability theorems for incompressible flows. The main theorems of 
this article are stated in Section 4, and proved in Section 5.

\section{B\'enard Problem}

\subsection{Boussinesq equations}
The B\'enard problem can be modeled by the Boussinesq equations. 
In this paper, we consider the B\'enard problem in a two-dimensional (2D) 
domain $\mathbb R^1 \times (0, h) \subset \R^2$ ($h > 0$). 
The Boussinesq equations, which govern the motion and states of the 
fluid flow, are as follows; 
see among others Rayleigh \cite{rayleigh}, Drazin and Reid \cite{dr}  and 
Chandrasekhar \cite{chandrasekhar}:
\begin{align}
\label{eq2.1}
& \frac{\partial  u}{\partial t} + (u\cdot \nabla )u - \nu\Delta u
+ \rho_0^{-1} \nabla p = -g k[1-\alpha(T-\bar T_0)], \\
\label{eq2.2}
& \frac{\partial T}{\partial t} + (u\cdot \nabla )T - \kappa\Delta T=0, \\
\label{eq2.3}
& \divv\,\, u = 0,
\end{align}
where $\nu,\kappa,\alpha,g$ are constants defined as in
(\ref{eq1.1}), $u=(u_1,u_2)$ the velocity field, $p$ the pressure 
function, $T$ the temperature function, 
$\bar T_0$ and $\bar T_1$ constants representing the lower and upper 
surface temperatures at $x_2 =0, h$, and $k=(0,  1)$ 
the unit vector in $x_3$-direction.

To make the equations non-dimensional, let
\begin{eqnarray*}
x&=& hx', \\
t&=& h^2 t'/\kappa,\\
u&=& \kappa u'/h, \\
T&=& \beta h (T'/\sR) + \bar T_0 - \beta h x_2', \\
p &=& \rho_0 \kappa^2  p'/h^2  +p_0 
    -g \rho_0 ( h x_2' + \alpha \beta h^2 (x_2')^2/2), \\
\pr &=& \nu /\kappa.
\end{eqnarray*}
Here the Rayleigh number $R$ is defined by (\ref{eq1.1}),
and $\pr=\nu /\kappa$  is the Prandtl number.

Omitting the primes, the equations (\ref{eq2.1})-(\ref{eq2.3})
can be rewritten as follows
\begin{align}
\label{eq2.4}
& \frac{1}{\pr} 
   \left[\frac{\partial u}{\partial t} + (u \cdot \nabla ) u + \nabla p\right]
 -  \Delta u -  \sR T k =0, \\
\label{eq2.5}
& \frac{\partial T}{\partial t} + (u\cdot \nabla )T - \sR u_3 - \Delta T =0, \\
\label{eq2.6}
& \divv\,\, u = 0.
\end{align}

The non-dimensional domain is $\Omega =
\mathbb R^1 \times (0,1) \subset \R^2$.
We consider periodic boundary condition in the $x_1$-direction 
\be\label{eq2.7}
(u, T)(x_1, x_2, t) = (u, T)(x_1 + kL, x_2, t) \quad \forall k \in \mathbb Z.
\ee

At the top and bottom boundary ($x_2=0, 1$), 
different combinations of top and bottom boundary conditions are normally 
used in different physical setting such as {\it rigid-rigid}, 
{\it rigid-free}, {\it free-rigid}, and {\it free-free}.
For instance, we have 

{\sc Dirichlet boundary condition (rigid-rigid):}
\be 
\label{eq2.8}
T=0, \quad u=0 \quad \text{ at } x_2=0, 1.
\ee

{\sc Free-free boundary condition:}
\be 
\label{eq2.9}
T=0, \quad u_2=0 \quad \frac{\p u_1}{\p x_2} =0 \quad \text{ at } x_2=0, 1.
\ee

{\sc Free-rigid boundary condition:}
\be 
\label{eq2.10}
\left\{
\begin{aligned}
& T=0, \quad u=0 && \text{ at } x_2=0,\\
& T=0, \quad u_2 =0, \quad \frac{\p u_1}{\p x_2} =0 
&& \text{ at } x_2= 1.
\end{aligned}
\right.
\ee

The initial value conditions are given by  
\begin{equation}
\label{eq2.11}
(u, T) = (u_0, T_0) \qquad \text{ at } t=0.
\end{equation}

\subsection{Functional setting}

For simplicity, we proceed in this article with the  set of boundary 
conditions given by (\ref{eq2.7}) and (\ref{eq2.8}),
and similar results hold true as well for other combinations 
of boundary conditions.

Let
\begin{align*} 
& H=\{ (u, T) \in L^2(\Omega)^3  \ | \ 
\text{ div} u=0, u_2|_{x_2=0, 1}=0, u_1 \text{ is periodic 
in $x_1$-direction}\}, \\
& V= \{ (u, T) \in H^1_0(\Omega)^3  \ | \ 
\text{ div} u=0, (u,T) \text{ is periodic in $x_1$ 
direction}\},\\
& H_1 = V\cap H^2(\Omega)^3.
\end{align*}

Let $G :H_1 \to H,$  and 
$L_\lambda =-A + B_\lambda :H_1 \to H$  be defined by 
\begin{align*}
& G (\psi) = (- P[(u\cdot\nabla )u], 
                       - (u\cdot \nabla )T ),  \\
& A\psi = ( -P (\Delta u), -\Delta T ), \\
& B_\lambda \psi = \lambda ( P(Tk), u_2 ), 
\end{align*}
for any $\psi=(u, T) \in H_1$.  
Here $\lambda =\sR$, and $P$ the Leray projection to $L^2$ fields.

Then the Boussinesq equations (\ref{eq2.4})--(\ref{eq2.8})
can be rewritten in the following operator form
\begin{equation}\label{eq2.12}
\frac{d\psi}{dt} = L_\lambda \psi + G(\psi), \qquad \psi=(u,T).
\end{equation}

\subsection{Attractor bifurcation of the B\'enard problem}

Let $\{S_\lambda(t)\}_{t\ge 0}$ be an operator semi-group generated by
the equation (\ref{eq2.12}).
Then the solution of (\ref{eq2.12})  can be expressed
as 
\[
\psi(t, \psi_0) = S_\lambda(t)\psi_0, \qquad t\ge 0.
\]

\begin{definition}
\label{df2.1}
A set $\Sigma \subset H$ is called an invariant set of
(\ref{eq2.12}) if $S(t) \Sigma = \Sigma$ for any $t\ge 0$. An
invariant set $\Sigma \subset H$ of (\ref{eq2.12}) is called  an
attractor if $\Sigma$ is compact, and there exists a neighborhood $U
\subset H$ of $\Sigma$ such that for any $\psi_0\in U$ we have 
$$
\lim_{t\to \infty}\dist_H(\psi(t,\psi_0),\Sigma)= 0.
$$
\end{definition}

\begin{definition}
\label{df2.2}
\begin{enumerate}

\item We say that the equation (\ref{eq2.12}) bifurcates from $(\psi,\lambda) =
(0,\lambda_0)$ to invariant sets $\Omega_\lambda$, if there exists a
sequence of invariant sets $\{\Omega_{\lambda_n}\}$ of (\ref{eq2.12}) 
such that  $0 \notin \Omega_{\lambda_n}$ and 
\begin{align*}
& \lim_{n\to \infty} \lambda_n = \lambda_0,  \\
& \lim_{n\to \infty} \max_{x\in \Omega_{\lambda_n}} |x| =0. 
\end{align*}

\item If the invariant sets $\Omega_\lambda$ are attractors of
(\ref{eq2.12}), then the bifurcation is called attractor bifurcation.

\end{enumerate}
\end{definition}

We are now in position to state the attractor bifurcation theorem 
for the B\'enard problem (\ref{eq2.4})-(\ref{eq2.8}). 
The linearized equations of (\ref{eq2.4})-(\ref{eq2.6}) are given by 
\begin{equation}
\label{eq2.13}               
\left\{
\begin{aligned}
& - \Delta u + \nabla p - \sR Tk =0,  \\
& - \Delta T - \sR u_2 =0, \\
& \divv\,u =0, 
\end{aligned}\right.
\end{equation}
where $R$ is the Rayleigh number. These equations are supplemented with 
the same boundary conditions (\ref{eq2.7}) and (\ref{eq2.8}) 
as  the nonlinear Boussinesq system. 
This eigenvalue problem 
for the Rayleigh number $R$ is
symmetric. Hence, we know that all  eigenvalues ${R}_k$ with
multiplicities $m_k$ of (\ref{eq2.13}) with (\ref{eq2.7}) and (\ref{eq2.8}) 
are real numbers, and
$$
0< { R}_1 <\cdots <{ R}_k <{ R}_{k+1}<\cdots\,.
$$
The first eigenvalue $R_1$
is a function of 
the period $L$. The critical Rayleigh number $R_c$  is given by 
\be 
\label{eq2.14}
R_c = \min_{L>0} R_1(L).
\ee

Let the multiplicity
of $R_c$ be $m_1=m$ ($m=$ even), and the first eigenspace be denoted by 
$E_0$.
Then we have the following attractor bifurcation theorem.

\begin{theorem}\cite{benard,b-book}
\label{th2.3}
For the B\'enard problem (\ref{eq2.4}-\ref{eq2.8}), 
the following assertions hold true.

\begin{enumerate}
\item When $R \le R_c$, 
the steady state $(u,T) =0$ is a globally asymptotically stable
in $H$.

\item The equations bifurcate from $((u,T),R) = (0,R_c)$ to attractors 
$\Sigma_{R}$ for $R>R_c$, with $m-1\le \dim \Sigma_{R} \le m$, and $\Sigma_{R}$ is an ($m-1$) 
dimensional homological  sphere, i.e. $\Sigma_{R}$ has the same homology 
as $S^{m-1}$.

\item For any $(u, T) \in \Sigma_{R}$, the velocity field $u$ can be
expressed as 
$$
u = \sum^{m}_{k=1} \alpha_k e_k + o\left( \sum^{m}_{k=1}
\alpha_k e_k\right),
$$
where $e_k$ are the velocity fields of the first eigenvectors in $E_0$.

\item  For any  open bounded neighborhood $U \subset H$ of $(u, T)=0$,  
the attractor $\Sigma_{R}$ attracts $U\setminus \Gamma$  in $H$, 
where $\Gamma$ is the stable manifold of $(u, T)=0$ 
with co-dimension $m$ in $H$.

\end{enumerate} 
\end{theorem}

\begin{remark}
\label{rm2.4}
{\rm
Results similar to this attractor bifurcation 
theorem hold true as well for 
the 3D B\'enard problems; see \cite{benard, b-book}.
}
\end{remark}

\section{Preliminaries}
\subsection{Center manifold functions}
To study the structure of the bifurcated attractors of  
(\ref{eq2.4}-\ref{eq2.8}), it is necessary to consider 
the reduction of nonlinear evolution equations to 
center manifolds. 
To this end, we introduce in this section a method to derive a
first order approximation of the central manifold functions,
which was introduced and used in \cite{b-book}. 

Let  $H$   and  $H_1$ be two Hilbert spaces,
and let $H_1 \hookrightarrow H$ be a dense and compact inclusion.
We consider the following
nonlinear evolution equation
\be
\label{eq3.1}
\left\{
\begin{aligned}
& \frac{du}{dt} = L_\lambda u +G(u,\lambda), \\
& u(0) = u_0,
\end{aligned}
\right.
\ee
where $u: [0, \infty) \to H$  is the unknown function, $\lambda \in
\mathbb R$  is the  system  parameter, and
$L_\lambda:H_1\to H$ are parameterized linear completely
continuous fields depending continuously on $\lambda\in \R^1$, which
satisfy
\begin{equation}
\label{eq3.2}
\left\{\begin{aligned}
& L_\lambda = -A + B_\lambda && \text{ is a sectorial operator}, \\
& A:H_1 \to H && \text{ a linear homeomorphism}, \\
& B_\lambda :H_1\to H && \text{ parameterized linear compact
operators.}
\end{aligned}\right.
\end{equation}
It is easy to see \cite{henry, pazy} that $L_\lambda$
generates an analytic semi-group $\{e^{-tL_\lambda}\}_{t\ge 0}$.
Then we can define  fractional power operators $L^\alpha_\lambda$ for any
$0\le \alpha \le 1$ with domain $H_\alpha = D(L^\alpha_\lambda)$ such that
$H_{\alpha_1} \subset H_{\alpha_2}$ if $\alpha_1 > \alpha_2$, and $H_0=H$.

Furthermore, we assume that the nonlinear term
$G(\cdot, \lambda):H_\alpha \to H$, for some $0\le \alpha < 1$, 
is a family of parameterized $C^r$
bounded operators ($r\ge 1$) continuously depending on the parameter
$\lambda\in \R^1$, such that
\begin{equation}
\label{eq3.3}
 G(u,\lambda) = o(\|u\|_{H_\alpha}), \quad \forall\,\, \lambda\in \R^1.
\end{equation}

In this paper, we are interested in the case where 
$L_\lambda = -A +B_\lambda$ are sectorial operators such
 that there exist an eigenvalue sequence $\{\rho_k\}
\subset \C^1$ and an eigenvector sequence $\{e_k, h_k\}\subset
H_1$ of $A$:
\begin{equation}
\label{eq3.4}
\left\{\begin{aligned}
& Az_k = \rho_kz_k,  \qquad z_k=e_k + i h_k, \\
& \text{Re} \rho_k\to \infty \,\,(k\to\infty), \\
& |\text{Im} \rho_k / (a \text{Re} \rho_k) | \le c,
\end{aligned}\right.
\end{equation}
for some $a, c > 0$, such that
$\{e_k, h_k\}$ is a basis of $H$.

Condition (\ref{eq3.4})  implies that $A$ is a sectorial operator.
For the operator $B_\lambda:H_1\to H$, we also assume that
there is a constant $0<\theta<1$ such that
\begin{equation}
\label{eq3.5}
B_\lambda :H_\theta \longrightarrow H \,\,\text{bounded, $\forall$
$\lambda\in \R^1$.}
\end{equation}
Under conditions (\ref{eq3.4}) and (\ref{eq3.5}), the operator
$L_\lambda=-A + B_\lambda$ is a sectorial operator.

Let $H_1$  and $H$ be decomposed into
\begin{equation}
\label{eq3.6} \left\{
\begin{aligned}
& H_1 = E^\lambda_1 \oplus E^\lambda_2,  \\
& H = \widetilde E^\lambda_1 \oplus \widetilde E^\lambda_2,
\end{aligned}\right.
\end{equation}
for $\la$ near $\lambda_0 \in \R^1$, where $E^\lambda_1$,
$E^\lambda_2$ are invariant subspaces of $L_\lambda$, such that
\begin{align*}
&\dim E^\lambda_1<\infty, \\
& \tilde E^\lambda_1 = E^\lambda_1, \\
& \widetilde E^\lambda_2 = \text{closure of $E^\lambda_2$ in $H$.}
\end{align*}
In addition,  $L_\lambda$ can be decomposed into $L_\lambda =
\cL^\lambda_1 \oplus \cL^\lambda_2$ such that for any $\lambda$
near $\lambda_0$,
\begin{equation}
\label{eq3.7}
\begin{cases}
\cL^\lambda_1 = L_\lambda |_{E^\lambda_1} : E^\lambda_1
\longrightarrow \widetilde E^\lambda_1, & \\
\cL^\lambda_2 = L_\lambda|_{E^\lambda_2}:E^\lambda_2 \longrightarrow
\widetilde E^\lambda_2, &
\end{cases}
\end{equation}
where all eigenvalues of $\cL^\lambda_2$ possess negative real
parts, and the eigenvalues of $\cL^\lambda_1$ possess  nonnegative
real parts at $\lambda=\lambda_0$.

Thus, for $\lambda$ near $\lambda_0$, equation (\ref{eq3.1}) can be
written as
\begin{equation}
\label{eq3.8} 
\left\{
\begin{aligned}
& \frac{dx}{dt} = \cL^\lambda_1 x +G_1(x,y,\lambda), & \\
& \frac{dy}{dt} = \cL^\lambda_2 y + G_2(x,y,\lambda), &
\end{aligned}
\right.
\end{equation}
where $u=x+y \in H_1$, $x\in E^\lambda_1$, $y\in E^\lambda_2$,
$G_i(x,y,\lambda) = P_iG(u,\lambda)$, and $P_i:H\to \widetilde E_i^\la$
are canonical projections.
Furthermore, let
$$E^\la_2(\alpha)=\text{ closure of $E^\la_2$ in } H_\alpha,
$$
with $\alpha < 1$ given by (\ref{eq3.3}).

The following center manifold theorem is classical; see \cite{henry, temam}.

\begin{theorem}
\label{th3.1} 
Assume (\ref{eq3.3})--(\ref{eq3.7}) hold true. 
Then there exists a neighborhood of $\lambda_0$ given
by $|\lambda-\lambda_0|<\delta$ for some $\delta>0$, a
neighborhood $B_\lambda \in E^\lambda_1$ of $x=0$, and a $C^1$
function $\Phi(\cdot,\lambda): B_\lambda \to E^\lambda_2(\theta)$
depending continuously on $\lambda$, such that
\begin{enumerate}
\item $\Phi(0,\lambda) =0$, $\Phi'_x(0,\lambda) =0$,

\item the set
$$
M_\lambda = \left\{(x,y)\in H\mid x\in B_\lambda, \,\,y =
\Phi(x,\lambda) \in E^\lambda_2(\theta)\right\},
$$
called the center manifolds, are locally invariant for
(\ref{eq3.1}), i.e.\ for each $u_0 \in M_\lambda$
\[
u_\lambda(t,u_0) \in M_\lambda, \qquad \forall\,\, 0\le t<t(u_0)
\]
for some $t(u_0)>0$, where $u_\lambda(t,u_0)$ is the solution of
(\ref{eq3.1});

\item if $(x_\lambda(t),y_\lambda(t))$ is a solution of (\ref{eq3.8}),
then there is a $\beta_\lambda>0$ and $k_\lambda>0$ with $k_\lambda$
depending on $(x_\lambda(0), y_\lambda(0))$ such that
\[
\| y_\lambda(t) - \Phi(x_\la(t),\lambda) \|_H\le k_\lambda
e^{-\beta_\lambda t}.
\]
\end{enumerate}
\end{theorem}

Also, it is classical that to bifurcation problem of (\ref{eq3.1}) 
is reduced to that for the following finite dimensional 
system:
\be 
\label{eq3.9}
\frac{dx}{dt} = \cL^\lambda_1 x +g_1(x, \Phi_\la(x),\lambda), 
\ee
for $x \in B_\la \subset E_1^\la$.

Now we give a formula to calculate the center manifold functions. 
Let the nonlinear operator $G$ be given by
\be
\label{eq3.10}
G(u, \la)=G_k(u, \la) + o(|u|^k),
\ee
for $k\ge 2$, where $G_k(u, \la)$ is a $k$-multilinear operator:
\begin{align*}
&G_k: H_1 \times \cdots \times H_1 \to H, \\
& G_k(u, \la) = G_k(u,\cdots, u, \la).
\end{align*}

The following theorem was proved in \cite{b-book}.

\begin{theorem}
\label{th3.2} 
Under the conditions (\ref{eq3.3})-(\ref{eq3.7}) and (\ref{eq3.10}), 
the center manifold function $\Phi(x, \la)$ can be expressed as
\be\label{eq3.11} 
\Phi(x, \la) = (-\cL^\lambda_2)^{-1}P_2 G_k(x,
\la) + O(|\text{\rm Re}\beta(\la)| \cdot \|x\|^k) + o(\|x\|^k), \ee where
$\cL^\lambda_2$ is given by (\ref{eq3.7}), $P_2:H\to \tilde E_2$ the
canonical projection, $x \in E_1^\la$, and
$\beta(\la)=(\beta_1(\la), \cdots, \beta_m(\la))$ the eigenvalues
of $\cL^\lambda_1$.
\end{theorem}

\begin{remark}
\label{rm3.3}
{\rm
Consider the case where $L_\la: H_1 \to H$ is symmetric. 
Then the eigenvalues are real, and the eigenvectors form an orthogonal 
basis of $H$. Therefore, we have 
\begin{align*}
& u=x+y \in E_1^\la \oplus E_2^\la,\\
& x = \sum^m_{i=1} x_i e_i \in E_1^\la, \\
& y = \sum^\infty_{i=m+1} x_i e_i \in E_2^\la.
\end{align*}
Then near $\la=\la_0$, the formula 
(\ref{eq3.11}) can be expressed  as follows. 
\be 
\label{eq3.12}
\Phi(x, \la) = \sum^{\infty}_{j=m+1} \Phi_j(x, \la) e_j 
+ O(|\text{\rm Re}\beta(\la)| \cdot \|x\|^k) + o(\|x\|^k), \ee
where 
\begin{align*}
& \Phi_j(x, \la)= -\frac{1}{\beta_j(\la)} 
\sum_{1\le j_1, \cdots, j_k \le m} a_{j_1 \cdots j_k}^j 
x_{j_1}\cdots x_{j_k}, \\
& a_{j_1 \cdots j_k}^j = (G_k(e_{j_1}, \cdots,e_{j_k}, \la), e_j)_H. 
\end{align*}
In many applications, the coefficients $a_{j_1 \cdots j_k}^j$ can be computed, 
and the first $m$ eigenvalues $\beta_1(\la), \cdots, \beta_m(\la)$
satisfy 
$$|\text{\rm Re}\beta(\la_0)| =
\sqrt{\sum^m_{j=1}(\text{\rm Re}\beta_j(\la_0))^2}=0.
$$
Hence (\ref{eq3.12}) gives an explicit 
formula for the first approximation of the 
center manifold functions.
}
\end{remark}

\subsection{$S^1$-attractor bifurcation}
In this section, we study the structure of the bifurcated attractor 
of (\ref{eq3.9}) when $m=2$. Namely, 
we consider a two-dimensional system as follows:
\be
\label{eq3.13}
\frac{dx}{dt}=\beta(\la) x - g(x, \la), \quad x \in \mathbb R^2.
\ee
Here $\beta(\la)$  is a continuous function of $\la$ 
satisfying 
\be
\label{eq3.14}
\beta(\la)\left\{
\begin{aligned}
& < 0 && \quad \text{ if } \la < \la_0, \\
& = 0 && \quad \text{ if } \la = \la_0, \\
& > 0 && \quad \text{ if } \la > \la_0, 
\end{aligned}
\right.
\ee
and 
\be
\label{eq3.15}
\left\{
\begin{aligned}
& g(x, \la)=g_k(x, \la) + o(|x|^k), \\
& g_k(\cdot, \la) \text{ is a k-multilinear field}, \\
& C_1 |x|^{k+1} \le (g_k(x, \la), x),
\end{aligned}
\right.
\ee
for some integer $k=2m+1 \ge 3$, and some constants $0< C_1$.

The following theorem was proved in \cite{b-book}, which shows that 
under conditions (\ref{eq3.14})  and (\ref{eq3.15}), the system 
(\ref{eq3.13}) bifurcates to an $S^1$-attractor.
\begin{figure}
        \centering \includegraphics[height=.5\hsize]{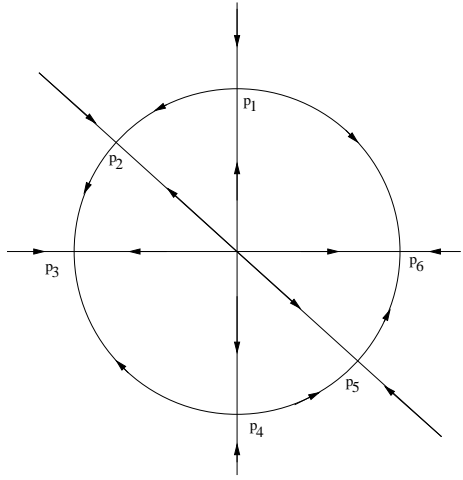}
\caption{$\Omega_\la$ has $4N + n$ ($N=1$ and $n=2$ shown here)
singular points, where $p_1$ and $p_4$ are saddles, 
$p_3$ and $p_6$ are nodes, and
$p_2$ and $p_5$ are singular points with index zero.}
\label{fg3.1}
\end{figure}
\begin{theorem}
\label{th3.4}
Let the condition (\ref{eq3.14})  and (\ref{eq3.15}) hold true. Then 
the system  (\ref{eq3.13}) 
bifurcates from $(x, \la)=(0,\la_0)$ to an attractor
$\Sigma_\la$, which is homeomorphic to $S^1$, 
for $\la_0 < \la< \la_0+\varepsilon$ for some $\varepsilon >0$. Moreover,
one and only one of the following is true.

\begin{enumerate}
\item $\Sigma_\la$ is a periodic orbit,

\item $\Sigma_\la$ consists of infinite number of singular points, or

\item $\Sigma_\la$ contains at most $2(k+1)=4(m+1)$ singular points, consisting of $2 N$  saddles points,  $2N$ stable node points, and $n$ ($\le 4(m+1) - 4N$) singular points with index zero, as shown in Figure~\ref{fg3.1}. 

\end{enumerate}
\end{theorem}

\subsection{Structural stability theorems}
In this subsection, we recall some results on structural stability 
for 2D divergence-free vector fields developed in \cite{amsbook}, 
which are 
crucial to study the asymptotic structure in the physical space 
of the bifurcated  solutions of the B\'enard problem.

Let $C^r(\Omega,\R^2)$ be the space of all $C^r$ $(r\ge 1)$ vector
fields on $\Omega=\mathbb R^1 \times (0, 1)$, 
which are periodic in $x_1$ direction with period $L$, 
let $D^r(\Omega,\R^2)$ be the space of all $C^r$  
divergence-free vector
fields on $\Omega=\mathbb R^1 \times (0, 1)$, 
which are periodic in $x_1$ direction with period $L$, and with no normal 
flow condition in $x_2$-direction:
\[
D^r(\Omega, \R^2) 
= \left\{ v\in C^r(\Omega,\R^2)\mid v_2 = 0 \text{ at } x_2=0, 1\ \right\}.
\]
Furthermore, we let 
\begin{align*}
& B^r_0(\Omega, \R^2) 
= \left\{ v\in D^r(\Omega,\R^2)\mid v = 0 \text{ at } x_2=0, 1\ \right\}, \\
& B^r_1(\Omega, \R^2) 
= \left\{ v\in D^r(\Omega,\R^2)\quad \Big|\quad  
\begin{matrix}
\displaystyle  v = 0 \text{ at } x_2=0 \\
\displaystyle  v_2=\frac{\partial v_1}{\partial x_2}=0 \text{ at } x_2=1
     \end{matrix} \ \right\}.
\end{align*}

\begin{definition}
\label{df3.5}
Two vector fields $u,v \in C^r(\Omega, \R^2)$ are called topologically
equivalent if there exists a homeomorphism of $\varphi:\Omega\to
\Omega$, which takes the orbits of $u$ to orbits of $v$ and preserves
their orientation.
\end{definition}

\begin{definition}
\label{df3.6}
Let $X= D^r(\Omega, \R^2)$  or $X= B^r_0(\Omega, \R^2)$. 
A vector field $v\in X$ is called structurally stable
in $X$ if there exists a neighborhood $U\subset X$ of $v$ 
such that for any $u\in U$, $u$ and $v$ are
topologically equivalent.
\end{definition}

Let $v\in D^r(\Omega, \R^2)$. 
We recall next some basic facts and definitions on divergence--free
vector fields. 
\begin{enumerate}
\item A point $p\in \Omega$ is called a singular point of $v$ if
$v(p) =0$; a singular point $p$ of $v$ is called non-degenerate if
the Jacobian matrix $Dv(p)$ is invertible; $v$ is called regular if
all singular points of $v$ are non-degenerate.

\item An interior non-degenerate singular point of $v$ can be either
a center or a saddle, and a non-degenerate boundary singularity must
be a saddle.

\item Saddles of $v$ must be connected to saddles. An interior saddle
$p\in \Omega$ is called self-connected if $p$ is connected only to
itself, i.e., $p$ occurs in a graph whose topological form is that of
the number 8.
\end{enumerate}

Let $v\in B_0^r(\Omega, \R^2)$; then we know that 
each point on $x_2=0, 1$ is a singular point of $v$ in the usual sense. 
To study the structure of $v$, we need to classify the  boundary points 
as follows.

\begin{definition}
Let $u \in B^r_0 (TM)(r \ge 2)$.

\begin{enumerate}
\item  A point $p \in \partial M$ is called a $\partial$-regular point of
$u$ if $\frac{\partial u_\TT (p)}{\partial n}\not= 0$; otherwise,
$p \in \partial M$ is called a $\partial$-singular point of $u$.

\item  A $\partial$-singular point $p \in \partial M$ of $u$ is
called non-degenerate if 
\begin{equation} \label{7.1}
\det 
\left( \begin{matrix}
\displaystyle \frac{\partial^2 u_\TT (p)}{\partial \tau \partial n}
&\displaystyle \frac{\partial^2 u_\TT (p)}{\partial n^2}\\
\\
\displaystyle \frac{\partial^2 u_n(p)}{\partial \tau \partial n}
&\displaystyle \frac{\partial^2 u_n(p)}{\partial n^2}
\end{matrix}
\right) \not= 0.
\end{equation}
A non-degenerate $\partial$-singular point of $u$ is also called a
$\partial$-saddle point of $u$.

\item $u \in B^r_0(TM)$  $(r \ge 2)$ is called
$D$-regular if    a) $u$ is regular in $\cM$, and
b) all $\partial$-singular points of $u$ on $\partial M$
are non-degenerate.
\end{enumerate}
\end{definition}

The following theorem provides necessary and sufficient conditions 
for structural stability of a
divergence--free vector field.

\begin{theorem} \label{th3.8} {\rm \cite{amsbook}}
Let $u\in B^r_0(TM)(r\ge 2)$.  Then $u$ is
structurally stable in $B^r_0(TM)$ if and only if
\begin{enumerate}
\item[1)]  $u$ is $D$-regular;

\item[2)]  all interior saddle points of $u$ are self-connection; and

\item[3)]  each $\P$-saddle point of $u$ on $\P M$ is connected to a
$\P$-saddle point on the same connected component of $\P M$.
\end{enumerate}
Moreover, the set of all structurally stable vector fields is open
and dense in $B^r_0(TM)$.
\end{theorem}

\begin{remark}\label{th3.9}
{\rm 
For vector fields with free-rigid boundary conditions, the conditions 
for structural stability differs slightly. More precisely, 
$u\in B^r_1(TM)(r\ge 2)$ is structurally stable in $B^r_1(TM)$ if and only if
\begin{enumerate}
\item[1)]  all singular of $u$ in $\Omega$ and on $x_2=1$ are regular, and 
all $\p$-singular points on $x_2=0$ are $\p$-regular;

\item[2)]  all interior saddle points of $u$ are self-connected; and

\item[3)]  each saddle of $u$ on $x_2=1$ is connected to saddles 
on $x_2=1$, and  each $\P$-saddle point of $u$ on $x_2=0$ is connected to a
$\P$-saddle point on $x_2=0$.
\end{enumerate}
}
\end{remark}

\begin{remark}\label{th3.10}
{\rm 
For vector fields satisfying free-free boundary conditions, we set
\begin{align*}
& B^r_2(\Omega, \R^2) 
= \left\{ v\in D^r(\Omega,\R^2)\ \Big| \  
v_2=\frac{\partial v_1}{\partial x_2}=0 \text{ at } x_2=0, 1\ \right\}, \\
& B^r_3(\Omega, \R^2) 
= \left\{ v\in D^r_2(\Omega,\R^2)\quad \Big|\quad \int_{\Omega} u dx =0\ 
\right\}.
\end{align*}
Then 
$u\in B^r_2(\Omega, \R^2)$ (resp.$u\in B^r_3(\Omega, \R^2)$) 
is structurally stable in $B^r_2(\Omega, \R^2)$ (resp. 
in $B^r_2(\Omega, \R^2)$) if and only if
\begin{enumerate}
\item[1)]  $u$ is regular;

\item[2)]  all interior saddle points of $u$ are self-connected; and

\item[3)]  each boundary saddle of $u$ is connected to boundary saddles 
on the same connected component of $\P \Omega$ (resp. 
each boundary saddle of $u$ is connected to boundary saddles
not necessarily on the same connected component).
\end{enumerate}
The difference between these two cases is due to the zero-average condition
in the definition in $B^r_3(\Omega, \R^2)$, 
which implies that $B^r_3(\Omega, \R^2)$ does not contain the 
 harmonic field $v_0=(\alpha, 0)$ for any constant $\alpha\not=0$. Hence, 
an orbit connecting two saddles on different components of the boundary 
can not be broken with a perturbation in  $B^r_3(\Omega, \R^2)$ into orbits 
connecting only saddles on the same connected component of the boundary.
}
\end{remark}

\section{Structure of Bifurcated Solutions for the B\'enard Problem}
In this section, we study the topological structure of the bifurcated 
attractor and the asymptotic structure of solutions for the B\'enard problem.
It is known that for each type of boundary conditions, there is a minimal 
period $L_c$ satisfying (\ref{eq2.14}). Hereafter, we always take $L_c$ to 
be the period of (\ref{eq2.7}).

The main theorem in this article is as follows.
\begin{figure}
        \centering \includegraphics[height=.6\hsize]{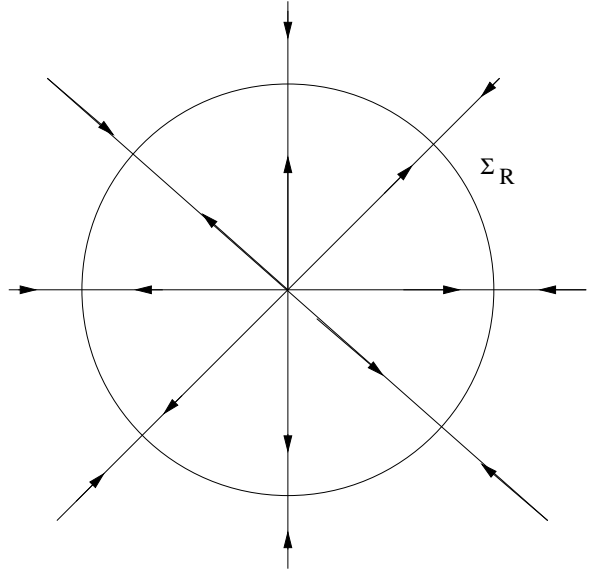}
\caption{All points on the bifurcated attractor $\Sigma_R=S^1$ are steady 
state solutions.}
\label{fg4.1}
\end{figure}

\begin{figure}
        \centering \includegraphics[height=.7\hsize]{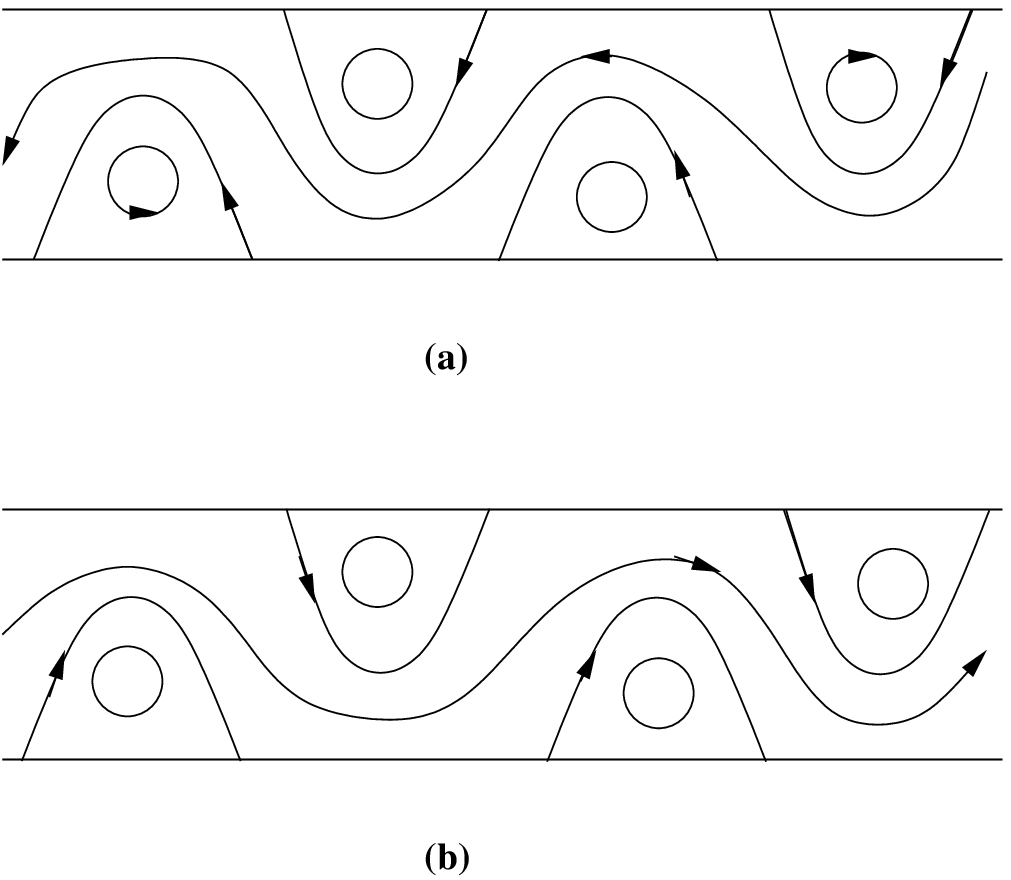}
\caption{Here the horizontal axis is the $x_1$-axis, and the vertical axis is 
the $x_2$-axis. With the Dirichlet boundary conditions on $x_2=0, 1$, the flow is not moving on both the top $x_2=1$  and the bottom $x_2=0$ boundaries.}
\label{fg4.2}
\end{figure}

\begin{theorem}
\label{th4.1}
For the B\'enard problem (\ref{eq2.4})-(\ref{eq2.8}), 
the following assertions hold true.

\begin{enumerate}

\item For $R> R_c$, the equations bifurcate from the trivial solution 
$((u, T), R) = (0, R_c)$ to an attractor $\Sigma_R$, homeomorphic to 
$S^1$, which consists of steady state solutions as shown in 
Figure~\ref{fg4.1}, where $R_c$  is the critical Rayleigh number.

\item For any $\psi_0=(u_0, T_0) \in H\setminus (\Gamma \cup E) $, 
there exists a time $t_0 \ge 0$ such that for any $t\ge t_0$, 
the vector field $u(t, \psi_0)$ is topologically 
equivalent to the structure as shown in either Figure~\ref{fg4.2}(a) or 
Figure~\ref{fg4.2}(b), where $\psi=(u(t, \psi_0), T(t, \psi_0))$  is 
the solution 
$\psi=(u(t, \psi_0), T(t, \psi_0))$  of (\ref{eq2.4})-(\ref{eq2.8}) 
with initial data $\psi_0$, $\Gamma$ is the stable manifold of 
the trivial solution $(u, T)=0$ with co-dimension $2$, and 
$$E=\left\{ (u, T) \in H\ | \ \int^1_0 u_1 dx_2 =0\right\}.$$

\end{enumerate}
\end{theorem}

The zonally moving meandering flow shown in Figure~\ref{fg4.2} 
appears often in many physical problems such as the Bransdator-Kushnir waves
in atmospheric circulation \cite{branstator, kushnir}.

\medskip

Theorem~\ref{th4.1} is also valid for the B\'enard problem 
(\ref{eq2.4})-(\ref{eq2.7}) with (\ref{eq2.10}). However the case with 
the free-free boundary condition is different.
More precisely, for the free-free boundary condition, it is easy to see that 
for any constant $\alpha$, the harmonic field 
$\psi_0=((\alpha, 0), 0)$ is a  solution of 
(\ref{eq2.4})-(\ref{eq2.6}). Therefore, we have to consider the problem  
(\ref{eq2.4})-(\ref{eq2.7}) with (\ref{eq2.9}) in the following 
function spaces:
\begin{align*} 
& \tilde H=\{ (u, T) \in L^2(\Omega)^3  \quad | \quad 
\text{ div} u=0, \int_\Omega u dx =0 \}, \\
& \tilde H_1 = \{ (u, T) \in \tilde H \cap H^2(\Omega)^3 \text{ satisfies 
(\ref{eq2.7}) and (\ref{eq2.9}) } \}
\end{align*}

Then we have the following theorem.

\begin{figure}
        \centering \includegraphics[height=.2\hsize]{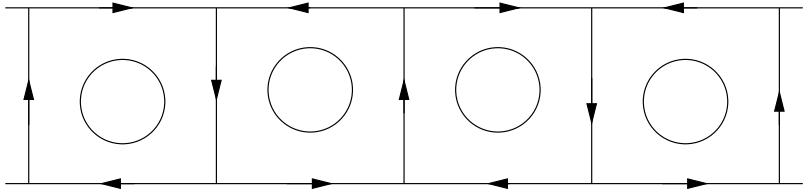}
\caption{Here  the horizontal axis is the $x_1$-axis, and the vertical axis is 
the $x_2$-axis. With the free slip boundary conditions on $x_2=0, 1$, the flow does move on both the top $x_2=1$  and the bottom $x_2=0$ boundaries.}
\label{fg4.3}
\end{figure}

\begin{theorem}
\label{th4.2}
For the B\'enard problem (\ref{eq2.4})-(\ref{eq2.7}) with 
boundary condition (\ref{eq2.9}), 
the following assertions hold true.

\begin{enumerate}

\item For $R> R_c$, the equations bifurcate from the trivial solution 
$((u, T), R) = (0, R_c)$ to an attractor $\Sigma_R$, homeomorphic to 
$S^1$, which consists of steady state solutions, 
where $R_c=27 \pi^4/4$  is the critical Rayleigh number.

\item For any $\psi_0=(u_0, T_0) \in \tilde H\setminus \Gamma $, 
there exists a time $t_0 \ge 0$ such that for any $t\ge t_0$, 
the vector field $u(t, \psi_0)$ is topologically 
equivalent to the structure as shown in Figure~\ref{fg4.3}, 
where $\psi=(u(t, \psi_0), T(t, \psi_0))$  is 
the solution 
$\psi=(u(t, \psi_0), T(t, \psi_0))$  of (\ref{eq2.4})-(\ref{eq2.7}) 
with (\ref{eq2.9}), $\Gamma$ is the stable manifold of 
the trivial solution $(u, T)=0$ with co-dimension $2$ in $\tilde H$.
\end{enumerate}

\end{theorem}

\section{Proof of Main Theorems}

\subsection{Eigenvectors of the linear Boussinesq equations}
We shall only prove Theorem~\ref{th4.1}. The proof of Theorem~\ref{th4.2}
is essentially the same, and we omit the details.
We proceed by first considering the eigenvalues and eigenvectors of 
the linearized equations of (\ref{eq2.4})--(\ref{eq2.6}):
\be
\label{eq5.1}
\left\{
\begin{aligned}
& \triangle u -\nabla p + \sqrt{R} T k = \beta(R)u, \\
&\triangle T + \sqrt{R} u_2 = \beta(R) T, \\
& \text{div} u =0,
\end{aligned}\right.
\ee
supplemented with the boundary conditions (\ref{eq2.7})  and (\ref{eq2.8}).

For $\psi=(u_1, u_2, T) \in H_1$, we take the separation of 
variables as follows:
\begin{align*}
& 
\psi= \left( 
-\sin\frac{2k\pi x_1}{L} h'(x_2), 
\frac{2k\pi }{L}\cos \frac{2k\pi x_1}{L} h(x_2),
\cos \frac{2k\pi x_1}{L} \theta(x_2)
\right), \\
& 
 \tilde \psi= \left( 
\cos \frac{2k\pi x_1}{L} h'(x_2),
\frac{2k\pi }{L}\sin \frac{2k\pi x_1}{L} h(x_2),
\sin \frac{2k\pi x_1}{L} \theta(x_2)
\right).
\end{align*}
Then it follows from (\ref{eq5.1}) that 
$(h, \theta)$ satisfies the following differential equations
\be
\label{eq5.2}
\left\{
\begin{aligned}
& \left(\frac{d^2}{dx_2^2}-a_k^2\right)^2 h - \sqrt{R} a_k \theta = \beta (R) 
\left(\frac{d^2}{dx_2^2}-a_k^2\right) h, \\
& - \left(\frac{d^2}{dx_2^2}-a_k^2\right) \theta - \sqrt{R} a_k h = - 
\beta (R) \theta, 
\end{aligned}\right.
\ee
supplemented with the following boundary conditions
\be
\label{eq5.3}\theta = 0, h=h'=0 \quad \text{ at } x_2=0, 1, 
\ee
where $a_k={2k\pi}/{L}$ and $L=L_c$ satisfies (\ref{eq2.14}).

The eigenvalue problem (\ref{eq5.2}) with (\ref{eq5.3}) is symmetric, 
and has a complete eigenvalue and eigenvector sequences for given $k$ and $R$:
\be
\label{eq5.4}
\left\{
\begin{aligned}
& \beta_{k1}(R) > \beta_{k2}(R) > \cdots, \\
& \lim_{j \to \infty}\beta_{kj}(R) = -\infty, \\
& h_{kj} \in H^4(0, 1)\cap H_0^2(0, 1)  \quad j=1, 2, \cdots, \\
& \theta_{kj} \in  H^2(0, 1)\cap H_0^1(0, 1) \quad j=1, 2, \cdots.
\end{aligned}\right.
\ee
Moreover, 
$$
\{(h_{kj},  \theta_{kj}) \ | \ j=1, 2, \cdots \}
$$
constitutes an orthogonal basis of $L^2(0, 1)\times L^2(0, 1)$.

Thus, we obtain the following 
complete set of eigenvectors for (\ref{eq5.1}) 
with boundary conditions (\ref{eq2.7})  and (\ref{eq2.8}):
\begin{align}
& 
\label{eq5.5}
\psi_{kj}= \left( 
-\sin\frac{2k\pi x_1}{L} h'_{kj}(x_2), 
\frac{2k\pi }{L}\cos \frac{2k\pi x_1}{L} h_{kj}(x_2),
\cos \frac{2k\pi x_1}{L} \theta_{kj}(x_2)
\right), \\
& \label{eq5.6}
 \tilde \psi_{kj}= \left( 
\cos \frac{2k\pi x_1}{L} h'_{kj}(x_2),
\frac{2k\pi }{L}\sin \frac{2k\pi x_1}{L} h_{kj}(x_2),
\sin \frac{2k\pi x_1}{L} \theta_{kj}(x_2)
\right).
\end{align}
where $0 \le k < \infty$, $1 \le j < \infty$. When  $k=0$, 
we derive from (\ref{eq5.1}), (\ref{eq5.5})  and (\ref{eq5.6})
that 
\be\label{eq5.7}
\left\{
\begin{aligned}
& \psi_{0j}= (0, 0, \sin j\pi x_2), \\
& \tilde \psi_{0j}=(h'_{0j}(x_2), 0, 0).
\end{aligned}\right.
\ee

\subsection{Singularity Cycle}
We shall show that the bifurcated attractor $\Sigma_R$ of 
(\ref{eq2.4})-(\ref{eq2.8}) given in Theorem~\ref{th2.3} contains a cycle of 
steady state solutions.

First, we note that the equations  
(\ref{eq2.4})-(\ref{eq2.8}) are invariant under the following translation:
$$\psi(x_1, x_2, t) \to \psi(x_1+ \alpha, x_2, t) \quad \forall \alpha 
\in \R.$$
Hence, if $\psi_0(x)$ is a steady state solution, then 
$\psi_0(x_1+ \alpha, x_2)$ are steady state solutions as well. 
It is easy to see that the set 
$$\mathcal S=\{ \psi_0(x_1 + \alpha, x_2)\ | \ \alpha 
\in \R\}$$
is a cycle $S^1$ in $H_1$. Therefore, each steady state of 
(\ref{eq2.4})-(\ref{eq2.8}) generates a cycle of steady state solutions.

Let 
\begin{align*}
& H'=\{ (u, T) \in H \ | \ u_1(-x_1, x_2)=-u_1(x_1, x_2)\}, \\
& H_1'=H_1\cap H'.
\end{align*}
It is easy to check that $H'$ and $H_1'$ are invariant spaces for 
the operator $L_\la + G$ given by (\ref{eq2.12}) in the sense that 
$$ L_\la + G: H_1'\to H',$$
where $ \la =\sqrt{R}$. On the other hand, it is clear that the sequence 
$\{\psi_{kj} \ | \ k=0, 1, \cdots, j=1, 2, \cdots \}$ defined by 
(\ref{eq5.5}) is  a basis of $H'$. Since the first eigenvalue $R_c$ of 
(\ref{eq2.13}) is simple, the first eigenvalue $\beta_1(R_c)$ of 
$L_\la$ in $H_1'$ is also simple, where
$$\beta_1(R_c)=\beta_{11}(R_c)=0,$$
and $\beta_{11}(R_c)$ is defined by (\ref{eq5.4}).
Hence by the classical Krasnoselskii bifurcation theorem, 
we know that the operator $L_\la + G$ bifurcates from 
$(\psi, \la)=(0, \sqrt{R_c})$ to a singular point in 
$H_1'$. Namely, the B\'enard problem (\ref{eq2.4})-(\ref{eq2.8})
bifurcates from $(\psi, R)=(0, R_c)$ a steady state 
solution. Therefore, the bifurcated attractor $\Sigma_R$  contains 
at least a cycle of steady state solutions.

\subsection{$S^1$-attractor: $\Sigma_R=S^1$}
To prove that $\Sigma_R=S^1$, by Theorem~\ref{th3.4}, it suffices 
to verify that the reduced equations of (\ref{eq2.4})-(\ref{eq2.8}) 
to the center manifold satisfy conditions 
(\ref{eq3.14})  and (\ref{eq3.15}).

For any $\psi=(u, T) \in H$, we have 
$$\psi = \sum_{k\ge 0, j \ge 1}\left(x_{kj}\psi_{kj} + 
y_{kj} \tilde \psi_{kj}    \right).$$
Since $L$ is the minimal period satisfying (\ref{eq2.14}), 
$\psi_{11}$  and $\tilde \psi_{11}$ are the first eigenvectors of 
(\ref{eq5.1}). Therefore the reduced equations of 
(\ref{eq2.4})-(\ref{eq2.8}) are given by 
\be\label{eq5.8}
\left\{
\begin{aligned}
& \frac{dx_{11}}{dt} = \beta_1(R) x_{11} + \frac{1}{\|\psi_{11}\|^2_H} 
(G(\psi, \psi), \psi_{11}), \\
& \frac{dy_{11}}{dt} = \beta_1(R) y_{11} + \frac{1}{\|\psi_{11}\|^2_H} 
(G(\psi, \psi),\tilde \psi_{11}).
\end{aligned}
\right.
\ee
Here for $\psi_1 =(u, T_1)$, $\psi_2 =(v, T_2)$, and $\psi_3=(w, T_3)$, 
$$
(G(\psi_1, \psi_2), \psi_{3})
=  - \int_{\Omega}\left[ (u\cdot \nabla v) w 
+(u\cdot \nabla T_2)T_3\right] dx.$$
Let the center manifold function be denoted by 
\be \label{eq5.9}
\Phi=\sum_{(k,j)\not= (1,1)}\left(
\Phi_{kj}(x_{11}, y_{11}) \psi_{kj} + 
\tilde \Phi_{kj}(x_{11}, y_{11}) \tilde \psi_{kj}    \right).
\ee
Note that for any 
$\psi_i \in H_1$ ($i=1, 2, 3$), 
\begin{align*}
& (G(\psi_1, \psi_2), \psi_{2}) =0, \\
& (G(\psi_1, \psi_2), \psi_{3}) =  - (G(\psi_1, \psi_3,), \psi_{2}).
\end{align*}
Then by $\psi=x_{11} \psi_{11} + y_{11}\tilde \psi_{11} + \Phi$, we have 
\begin{align}
(G(\psi, \psi), \psi_{11}) = 
  & (G(\tilde \psi_{11}, \tilde \psi_{11}), \psi_{11}) y_{11}^2 
      \label{eq5.10} \\
  &  - (G(\psi_{11}, \psi_{11}), \tilde \psi_{11}) x_{11} y_{11} 
           \nonumber \\
  &  - (G(\psi_{11}, \psi_{11}), \Phi)  x_{11} 
           \nonumber\\
  & + (G(\tilde \psi_{11}, \Phi)+ G(\Phi, \tilde \psi_{11}), \psi_{11}) y_{11}
           \nonumber \\
  & + (G(\Phi,\Phi), \psi_{11}), 
           \nonumber\\
(G(\psi, \psi), \tilde \psi_{11}) = 
  & (G(\psi_{11}, \psi_{11}), \tilde\psi_{11}) x_{11}^2 
      \label{eq5.11} \\
  &  - (G(\tilde\psi_{11}, \tilde\psi_{11}), \psi_{11}) x_{11} y_{11} 
           \nonumber \\
  &  - (G(\tilde\psi_{11}, \tilde\psi_{11}), \Phi)  y_{11} 
           \nonumber\\
  & + (G(\psi_{11}, \Phi)+ G(\Phi, \psi_{11}), \tilde\psi_{11}) x_{11}
           \nonumber \\
  & + (G(\Phi,\Phi),\tilde \psi_{11}). 
           \nonumber
\end{align} 
It is easy to check that for $k\not=0, 2$, 
\be\label{eq5.12}
\left\{
\begin{aligned}
& (G(\psi_{11}, \psi_{11}),\psi_{kj}) =0, \\
& (G(\tilde\psi_{11}, \tilde\psi_{11}),\psi_{kj}) =0, \\
& (G(\tilde\psi_{kj}, \tilde\psi_{11}),\psi_{11}) =0, \\
& (G(\tilde\psi_{kj}, \psi_{11}),\tilde \psi_{11}) =0,
\end{aligned}\right.
\ee
and for any $k\ge 0$, $j\ge 1$, we have
\be\label{eq5.13}
\left\{
\begin{aligned}
& (G(\psi_{11}, \psi_{11}),\tilde \psi_{kj}) =0, \\
& (G(\tilde\psi_{11}, \tilde\psi_{11}),\psi_{kj}) =0, \\
& (G(\psi_{kj}, \tilde\psi_{11}),\psi_{11}) =0, \\
& (G(\psi_{kj}, \psi_{11}),\tilde \psi_{11}) =0.
\end{aligned}\right.
\ee
By  (\ref{eq5.9}),  (\ref{eq5.12}) and (\ref{eq5.13}), the equalities 
(\ref{eq5.10}) and (\ref{eq5.11}) can be rewritten as 
\begin{align}
(G(\psi, \psi), \psi_{11}) = 
  & - \sum_{j=1}^\infty 
    [ (G(\psi_{11}, \psi_{11}), \psi_{0j}) \Phi_{0j} 
        + (G(\psi_{11}, \psi_{11}), \psi_{2j}) \Phi_{2j}] x_{11} 
      \label{eq5.14} \\
  &  - \sum_{j=1}^\infty
    [ (G(\tilde \psi_{11}, \psi_{11}), \tilde \psi_{0j})
         + (G(\tilde \psi_{0j}, \psi_{11}),  \psi_{11}) 
         y_{11} \tilde \Phi_{0j}
           \nonumber \\
  &  - \sum_{j=1}^\infty
    [ (G(\tilde \psi_{11}, \psi_{11}), \tilde \psi_{2j})
         + (G(\tilde \psi_{2j}, \psi_{11}), \tilde \psi_{11}) 
         y_{11} \tilde \Phi_{2j}
           \nonumber \\
  & + (G(\Phi,\Phi), \psi_{11}), 
           \nonumber
\end{align}
\begin{align}
(G(\psi, \psi), \tilde \psi_{11}) = 
  & - \sum_{j=1}^\infty 
    [ (G(\tilde \psi_{11}, \tilde\psi_{11}), \psi_{0j}) \Phi_{0j} 
        + (G(\tilde\psi_{11}, \tilde\psi_{11}), \psi_{2j}) \Phi_{2j}] y_{11} 
      \label{eq5.15} \\
  &  - \sum_{j=1}^\infty
    [ (G(\psi_{11}, \tilde\psi_{11}), \tilde \psi_{0j})
         + (G(\tilde \psi_{0j}, \tilde\psi_{11}), \tilde \psi_{11}) 
         x_{11} \tilde \Phi_{0j}
           \nonumber \\
  &  - \sum_{j=1}^\infty
    [ (G(\psi_{11}, \tilde\psi_{11}), \tilde \psi_{2j})
         + (G(\tilde \psi_{2j}, \tilde\psi_{11}), \psi_{11}) 
         x_{11} \tilde \Phi_{2j}
           \nonumber \\
  & + (G(\Phi,\Phi), \tilde\psi_{11}).
           \nonumber
\end{align}

Since the center manifold functions contains only higher order terms
$$\Phi(x_{11}, y_{11}) = O(|x_{11}|^2, |y_{11}|^2), $$
we derive that 
\be\label{eq5.16}
\left\{
\begin{aligned}
& (G(\Phi,\Phi), \psi_{11}) = o(|x_{11}|^3, |y_{11}|^3),\\
& (G(\Phi,\Phi), \tilde\psi_{11})= o(|x_{11}|^3, |y_{11}|^3).
\end{aligned}
\right.
\ee
Then direct calculation yields that
\begin{align}
\label{eq5.17}
(G(\tilde \psi_{2j}, \psi_{11}), \tilde \psi_{11}) = 
  & - (G(\tilde \psi_{2j}, \tilde \psi_{11}), \psi_{11}) \\
 = & \frac{\pi}{2} \int^1_0[ -h_{11}' (h_{2j}'h_{11}' + 
          2 h_{2j}h_{11}'')  \nonumber \\
    & + h_{11}  (h_{2j}'h_{11}   + 2 h_{2j}h_{11}')\nonumber  \\
    & + \theta_{11}  (h_{2j}'\theta_{11}   + 2 h_{2j}\theta_{11}')]dx_2
\nonumber \\
= & \frac{\pi}{2} \int^1_0\frac{d}{dx_2}(-h_{2j} (h_{11}')^2 + h_{2j}h_{11}^2 
       + h_{2j}\theta_{11}^2)dx_2\nonumber \\
=& 0.\nonumber 
\end{align}
It is clear that for any $k\ge 0$ and $j\ge 1$, 
\be
\|\psi_{kj}\|^2_H =\| \tilde \psi_{kj}\|^2_H.\label{eq5.18} 
\ee
Hence the reduced equations (\ref{eq5.8}) can be expressed as follows:
\begin{align}
\label{eq5.19}
\frac{dx_{11}}{dt} = 
& \beta_1(R) x_{11}- \frac{1}{\|\psi_{11}\|^2_H}
\sum_{j=1}^\infty 
    [ (G(\psi_{11}, \psi_{11}), \psi_{0j}) \Phi_{0j} x_{11} \\
& + (G(\psi_{11}, \psi_{11}), \psi_{2j}) \Phi_{2j} x_{11}
  + (G(\tilde \psi_{11}, \psi_{11}), \tilde \psi_{2j})
   y_{11} \tilde \Phi_{2j} \nonumber\\
& + (G(\tilde \psi_{11}, \psi_{11}), \tilde \psi_{0j}) y_{11} \tilde \Phi_{0j}
+ (G(\tilde \psi_{0j}, \psi_{11}), \tilde \psi_{11}) 
         y_{11} \tilde \Phi_{0j}]\nonumber \\
& + o(|x_{11}|^3, |y_{11}|^3), \nonumber\\
 \frac{dy_{11}}{dt} = 
& \label{eq5.20}
\beta_1(R) y_{11}- \frac{1}{\|\psi_{11}\|^2_H} 
\sum_{j=1}^\infty 
    [ (G(\tilde \psi_{11}, \tilde \psi_{11}), \psi_{0j}) \Phi_{0j} y_{11} \\
& + (G(\tilde \psi_{11}, \tilde \psi_{11}), \psi_{2j}) \Phi_{2j} y_{11}
  + (G(\psi_{11}, \tilde \psi_{11}), \tilde \psi_{2j})
   x_{11} \tilde \Phi_{2j} \nonumber\\
& + (G(\psi_{11}, \tilde \psi_{11}), \tilde \psi_{0j}) x_{11} \tilde \Phi_{0j}
+ (G(\tilde \psi_{0j}, \tilde \psi_{11}), \psi_{11}) 
         x_{11} \tilde \Phi_{0j}]\nonumber \\
& + o(|x_{11}|^3, |y_{11}|^3). \nonumber
\end{align}

By Theorem~\ref{th3.2} and (\ref{eq3.12}), we have 
\begin{align*}
\Phi_{0j} = & \frac{-1}{\|\psi_{0j}\|^2_H \beta_{0j}}
[ (G(\psi_{11}, \psi_{11}), \psi_{0j})  x^2_{11} 
+ (G(\tilde \psi_{11}, \tilde\psi_{11}), \psi_{0j})  y^2_{11}]\\
& + o(x_{11}^2 + y_{11}^2) + O(\beta_1(R) (x_{11}^2 + y_{11}^2)), \\
\Phi_{2j} = & \frac{-1}{\|\psi_{2j}\|^2_H \beta_{2j}}
[ (G(\psi_{11}, \psi_{11}), \psi_{2j})  x^2_{11} 
+ (G(\tilde \psi_{11}, \tilde\psi_{11}), \psi_{2j})  y^2_{11}]\\
& + o(x_{11}^2 + y_{11}^2) + O(\beta_1(R) (x_{11}^2 + y_{11}^2)), \\
\tilde \Phi_{0j} = & \frac{-1}{\|\tilde\psi_{0j}\|^2_H \beta_{0j}}
[ (G(\psi_{11}, \tilde\psi_{11}) + (G(\tilde \psi_{11}, \psi_{11}), 
     \tilde \psi_{0j}) x_{11} y_{11}]\\
& + o(x_{11}^2 + y_{11}^2) + O(\beta_1(R) (x_{11}^2 + y_{11}^2)), \\
\tilde\Phi_{2j} = & \frac{-1}{\|\tilde\psi_{2j}\|^2_H \beta_{2j}}
[ (G(\psi_{11}, \tilde\psi_{11}) + (G(\tilde \psi_{11}, \psi_{11}), 
  \tilde \psi_{2j})  x_{11} y_{11}]\\
& + o(x_{11}^2 + y_{11}^2) + O(\beta_1(R) (x_{11}^2 + y_{11}^2)), \\
\Phi_{kj} = &  o(x_{11}^2 + y_{11}^2) 
\qquad \forall k\not=0,2,\\
\tilde\Phi_{kj} = &  o(x_{11}^2 + y_{11}^2) 
\qquad \forall k\not=0,2,
\end{align*}
where $\beta_{kj}(R)$ are as in (\ref{eq5.4}). 

Also by direct computation, we obtain that 
 \begin{align*}
& (G(\psi_{11}, \psi_{11}), \psi_{0j})= 
(G(\tilde \psi_{11}, \tilde \psi_{11}), \psi_{0j}), \\
& (G(\psi_{11}, \psi_{11}), \psi_{2j})= 
-(G(\tilde \psi_{11}, \tilde \psi_{11}), \psi_{2j}), \\
& (G(\psi_{11}, \tilde \psi_{11})+ G(\tilde \psi_{11}, \psi_{11}),
\tilde \psi_{0j}) =0, \\
& (G(\psi_{11}, \tilde \psi_{11}),\tilde  \psi_{2j})
= (G( \tilde\psi_{11}, \psi_{11}),\tilde  \psi_{2j})
=  (G(\psi_{11}, \psi_{11}),  \psi_{2j}).
\end{align*}
Hence we have 
\begin{align}
\label{eq5.21} 
\Phi_{0j} 
   = & \frac{-1}{\|\psi_{0j}\|^2_H \beta_{0j}}
     (G(\psi_{11}, \psi_{11}), \psi_{0j}) ( x^2_{11}+ y^2_{11}) \\
     &    + o(x_{11}^2 + y_{11}^2) + O(\beta_1(R) (x_{11}^2 + y_{11}^2)), 
        \nonumber \\
\label{eq5.22}
\Phi_{2j} = & \frac{-1}{\|\psi_{2j}\|^2_H \beta_{2j}}
(G(\psi_{11}, \psi_{11}), \psi_{2j}) ( x^2_{11}- y^2_{11}) \\
& + o(x_{11}^2 + y_{11}^2) + O(\beta_1(R) (x_{11}^2 + y_{11}^2)), 
 \nonumber \\
\label{eq5.23}
\tilde \Phi_{0j} = & 
o(x_{11}^2 + y_{11}^2) + O(\beta_1(R) (x_{11}^2 + y_{11}^2)), \\
\label{eq5.24}
\tilde\Phi_{2j} = & \frac{-2}{\|\psi_{2j}\|^2_H \beta_{2j}}
(G(\psi_{11}, \psi_{11}), \psi_{2j})  x_{11} y_{11}\\
& + o(x_{11}^2 + y_{11}^2) + O(\beta_1(R) (x_{11}^2 + y_{11}^2)).\nonumber
\end{align}
Inserting (\ref{eq5.21})-(\ref{eq5.24}) into 
(\ref{eq5.19})  and (\ref{eq5.20}), we have 
\begin{align}
\label{eq5.25}
\frac{dx_{11}}{dt} = 
& \beta_1(R) x_{11}- \alpha x_{11}(x_{11}^2 + y_{11}^2)\\
& + o(x_{11}^3 + y_{11}^3) + O(\beta_1(R) (x_{11}^3 + y_{11}^3)),\nonumber\\
 \frac{dy_{11}}{dt} = 
& \label{eq5.26}
\beta_1(R) y_{11}- \alpha y_{11}(x_{11}^2 + y_{11}^2)\\
& + o(x_{11}^3 + y_{11}^3) + O(\beta_1(R) (x_{11}^3 + y_{11}^3)),\nonumber
\end{align}
where 
$$\alpha = 
\frac{-1}{\|\psi_{11}\|^2_H} 
\sum^\infty_{j=1} \left[
\frac{(G(\psi_{11}, \psi_{11}), \psi_{0j})^2}{\|\psi_{0j}\|^2_H \beta_{0j}(R)}
+ 
\frac{(G(\psi_{11}, \psi_{11}), \psi_{2j})^2}{\|\psi_{2j}\|^2_H \beta_{2j}(R)}
\right].
$$
We know that 
\begin{align*}
& \beta_1(R)=\beta_{11}(R) > \beta_{kj}(R) \quad \forall (k, j)\not=(1,1), \\
& \beta_{11}(R_c) =0.
\end{align*}
Hence near $R=R_c$, 
 $$\beta_{kj}(R) < 0 \quad \forall (k, j)\not=(1,1).$$
Consequently, $\alpha > 0$ and (\ref{eq5.25})  and (\ref{eq5.26})
satisfy (\ref{eq3.15}).

Also, we know that \cite{benard, b-book}
$$
\beta_1(R) \left\{\begin{aligned}
& < 0 \qquad && \text{if } R < R_c, \\
& = 0 \qquad && \text{if } R = R_c, \\
& > 0 \qquad && \text{if } R > R_c. 
\end{aligned}\right.
$$
Therefore (\ref{eq3.15}) holds true.

\subsection{Asymptotic structure of solutions}
By \cite{fmt}, we know that for any initial value $\psi_0 =(u_0, T_0) \in H$, 
there is a time $\tau \ge 0$ such that 
the solution $\psi =(u(t, \psi_0), T(t, \psi_0))$   is $C^\infty$ 
for $t> \tau$, and is uniformly bounded in $C^r$-norm for any given $r\ge 1$. 
Hence, by Theorem~\ref{th2.3}, we have 
\be
\label{eq5.27}
\lim_{t \to \infty} \min_{\phi \in \Sigma_R} \| \psi(t, \psi_0) - \phi\|_{C^r} 
=0.
\ee
We infer then from (\ref{eq5.25})  and (\ref{eq5.26}) that 
for any steady state solution $\Phi=(e, T) \in \Sigma_R$ of 
(\ref{eq2.4})-(\ref{eq2.8}), the vector field $e=(e_1, e_2)$ can be expressed 
as 
\be
\label{eq5.28}
\left\{
\begin{aligned}
& e_1 = r \cos\frac{2\pi}{L}(x_1 + \theta) \ h_{11}'(x_2) + 
    v_1(x_{11}, y_{11},\beta_1), \\
& e_2 = \frac{2\pi}{L} r \sin\frac{2\pi}{L}(x_1 + \theta) \ h_{11}(x_2) + 
    v_2(x_{11}, y_{11}, \beta_1), 
\end{aligned}\right.
\ee
for some $0\le \theta \le 2\pi$. Here 
\be
\label{eq5.29}
\left\{
\begin{aligned}
& r=\sqrt{x_{11}^2 +y_{11}^2 } = \sqrt{\beta_1(R)} + o(\sqrt{\beta_1(R)}) && 
 \text{ if } R > R_c, \\
& v_i(x_{11}, y_{11}, \beta_1) = o(\sqrt{\beta_1(R)}) && \text{ for }i=1, 2.
\end{aligned}\right.
\ee

\begin{figure}
        \centering \includegraphics[height=.5\hsize]{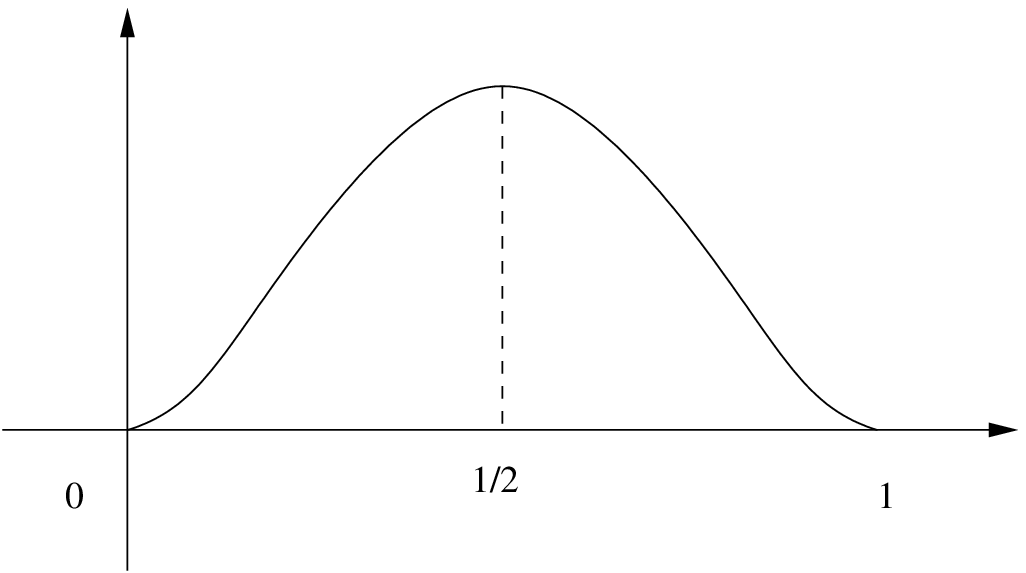}
\caption{}
\label{fg5.1}
\end{figure}

On the other hand, it is known that the first eigenfunction $h_{11}(x_2)$ of 
(\ref{eq5.2}) and (\ref{eq5.3}) at $R=R_c$ is given by 
\begin{align}
\label{eq5.30}
h_{11}(x_2) \simeq & \cos \alpha_0(x_2 -\frac{1}{2}) -0.06 
\cosh \alpha_1(x_2 -\frac{1}{2})\cos \alpha_2(x_2 -\frac{1}{2})\\
& + 0.1 \sinh \alpha_1(x_2 -\frac{1}{2})\sin \alpha_2(x_2 -\frac{1}{2}),
\nonumber
\end{align}
where $\alpha_0\simeq 3.97$, $\alpha_1\simeq 5.2$  and $\alpha_2\simeq 2.1$.
This function is schematically given by Figure~\ref{fg5.1}; see 
\cite{chandrasekhar}.

Now we show that the vector field 
\be
\label{eq5.31}
e_0 = \left(r \cos \frac{2\pi x_1}{L} h_{11}'(x_2),  
\frac{2\pi r}{L}\sin \frac{2\pi x_1}{L} h_{11}(x_2)  \right)
\ee
is $D$-regular in $\Omega=\mathbb R^1 \times (0, 1)$.

To this end, by (\ref{eq5.30}) we see that 
\begin{align*}
& h_{11}''(x_2)\not= 0 \qquad \text{ at } x_2=0, \frac12, 1, \\
&  h_{11}'(\frac12)=0,\qquad  h_{11}''(\frac12)\not=0.
\end{align*}
Hence 
$$ \text{det } De_0(x_1, x_2) \not=0, $$
for any $(x_1, x_2) =({kL}/{2}, 1/2)$  with 
$k=1, 2 , \cdots,$  and 
$$
 \text{det } 
\left(
\begin{matrix}
\displaystyle 
\frac{\partial^2 e^1_0}{\partial x_1 \partial x_2} & 
\displaystyle 
\frac{\partial^2 e^1_0}{\partial x_2^2}  \\
\\
\displaystyle 
\frac{\partial^2 e^2_0}{\partial x_1 \partial x_2} & 
\displaystyle 
\frac{\partial^2 e^2_0}{\partial x_2^2} 
\end{matrix}
\right)
=-\frac{2\pi r}{L} \sin^2\frac{2\pi x_1}{L} \ h_{11}''(x_2)  
\not=0, $$
for any $(x_1, x_2)=({(2k+1)L}/4, 0)$ or $(x_1, x_2)=({(2k+1)L}/4, 1)$
with $k=1, 2 , \cdots.$
Therefore, the vector field (\ref{eq5.31})  is $D$-regular, and consequently, 
the vector fields $e$ in (\ref{eq5.28}) are  $D$-regular  for any 
$R_c < R < R_c  + \varepsilon$ for some $\varepsilon > 0$ small.

Next we show that the following subspace of $H$
$$E=\{ \ (u, T)\in H_1 \ | \ \int_0^L\int^1_0 u_1 dx =0\}$$
is invariant for (\ref{eq2.4})-(\ref{eq2.8}). In fact, we can verify that 
\be
\label{eq5.32}
\int_0^L\int^1_0 P[(u\cdot \nabla)u] \cdot i dx=0 \quad \forall u \in H_1,
\ee
where $i=(0, 1)^t$ is the unit vector in the $x_1$-direction, and $P$ the 
Leray projection. Indeed, by the Helmholtz decomposition, we have 
$$[(u\cdot \nabla)u] \cdot i = P[(u\cdot \nabla)u] \cdot i 
+ \frac{\partial \phi}{\partial x_1},
$$
for some $\phi \in H^1(\Omega)$. Hence, 
$$\int_0^L\int^1_0 P[(u\cdot \nabla)u] \cdot i dx
= \int_0^L\int^1_0 (u\cdot \nabla)u_1 dx=0.$$

The invariance of $E$ for (\ref{eq2.4})-(\ref{eq2.8}) implies that 
for the vector field $e$ given in 
(\ref{eq5.28}), we have 
\be
\label{5.33}
\int^1_0 v_1 dx_2 =0.
\ee
Hence in the Fourier expansion of $e$ in (\ref{eq5.5}) and (\ref{eq5.6}), 
the coefficients of $\tilde \psi_{0j}$  are zero. By the connection lemma 
in \cite{amsbook}, it follows from (\ref{5.33}) that the vector field 
$e=(e_1, e_2)$  of (\ref{eq5.28}) is topologically 
equivalent to the vector field $e_0$ given by 
(\ref{eq5.31}), which has the topological structure as shown in 
Figure~\ref{fg4.3}.

For any initial value $\psi_0 =(u_0, T_0)\in H\setminus E$, 
\begin{align*}
& \psi_0 = \sum_k \alpha_k \tilde \psi_{0k} + \Phi_0, \\
& \Phi_0 \in E, \\
& \tilde \psi_{0k}=(\sin\pi x_2, 0, 0), \qquad k=2m+1, m=0, 1, \cdots.
\end{align*}
By (\ref{eq5.32}), we see that for the operator $G$ in (\ref{eq2.12}), 
$$(G(\psi), \tilde \psi_{0k})=0, 
$$
which implies that the solution of 
(\ref{eq2.4})-(\ref{eq2.8}) with (\ref{eq2.11}) has 
the following form
\be
\psi(t, \psi_0) = \sum_k \alpha_k e^{t\beta_{0k}}\tilde \psi_{0k} + 
\tilde \Phi(t, \psi_0), \label{eq5.34} 
\ee
where $\tilde \Phi \in E$, and $\beta_{0k} < 0$  near $R=R_c$.

Let 
$$K=\min\{ \ k \ | \ \alpha_k\not=0 \text{ and $k$ is odd}\},$$
and let $\psi(t, \psi_0) = ( u(t, \psi_0), T(t, \psi_0))$ be the solution 
of (\ref{eq2.4})-(\ref{eq2.8}) given by (\ref{eq5.34}). Then by 
(\ref{eq5.27}), the vector field $u(t, \psi_0)$ is topologically equivalent 
to the following vector field for any $t> t_0$ with $t_0>0$ sufficiently 
large
\be
\label{eq5.35}
\tilde u = e + (\alpha_k e^{-t\beta_{0k}}\sin k\pi x_2, 0),
\ee
where $e$ is as in (\ref{eq5.28}).

Since the vector field  $e$ is $D$-regular and topologically equivalent to 
the vector field as shown in Figure~\ref{fg4.3}. Then using the method 
for breaking saddle connections in \cite{amsbook}, it is easy to show that 
the vector field $\tilde u$ given by (\ref{eq5.35})  is topologically 
equivalent to the structure as shown in Figure~\ref{fg4.2}(a) 
if $\alpha_k <0$, and to the structure as shown in Figure~\ref{fg4.2}(b) 
if $\alpha_k >0$, for any $t>t_0$ sufficiently large. 

Thus the proof of Theorem~\ref{th4.1} is complete.


\begin{thebibliography}{10}

\bibitem{branstator}
{\sc G.~W. Branstator}, {\em A striking example of the atmosphere's leading
  traveling pattern}, J. Atmos. Sci., 44 (1987), pp.~2310--2323.

\bibitem{chandrasekhar}
{\sc S.~Chandrasekhar}, {\em Hydrodynamic and Hydromagnetic Stability}, Dover
  Publications, Inc.1981.

\bibitem{dr}
{\sc P.~Drazin and W.~Reid}, {\em Hydrodynamic Stability}, Cambridge University
  Press, 1981.

\bibitem{fmt}
{\sc C.~Foias, O.~Manley, and R.~Temam}, {\em Attractors for the {B}\'enard
  problem: existence and physical bounds on their fractal dimension}, Nonlinear
  Anal., 11 (1987), pp.~939--967.

\bibitem{henry}
{\sc D.~Henry}, {\em Geometric theory of semilinear parabolic equations},
  vol.~840 of Lecture Notes in Mathematics, Springer-Verlag, Berlin, 1981.

\bibitem{kirch}
{\sc K.~Kirchg{\"a}ssner}, {\em Bifurcation in nonlinear hydrodynamic
  stability}, SIAM Rev., 17 (1975), pp.~652--683.

\bibitem{kushnir}
{\sc K.~Kushnir}, {\em Retrograding wintertime low-frequency disturbances over
  the north pacific ocean}, J. Atmos. Sci., 44 (1987), pp.~2727--2742.

\bibitem{benard}
{\sc T.~Ma and S.~Wang}, {\em Dynamic bifurcation and stability in the
  {R}ayleigh-{B}\'enard convection}, Communication of Mathematical Sciences,
  2:2 (2004), pp.~159--183.

\bibitem{b-book}
\leavevmode\vrule height 2pt depth -1.6pt width 23pt, {\em Bifurcation Theory
  and Applications}, World Scientific, 2005.

\bibitem{amsbook}
\leavevmode\vrule height 2pt depth -1.6pt width 23pt, {\em Geometric Theory of
  Incompressible Flows with Applications to Fluid Dynamics}, vol.~119 of
  Mathematical Surveys and Monographs, American Mathematical Society,
  Providence, RI, 2005.

\bibitem{pazy}
{\sc A.~Pazy}, {\em Semigroups of linear operators and applications to partial
  differential equations}, vol.~44 of Applied Mathematical Sciences,
  Springer-Verlag, New York, 1983.

\bibitem{rabinowitz}
{\sc P.~H. Rabinowitz}, {\em Existence and nonuniqueness of rectangular
  solutions of the {B}\'enard problem}, Arch. Rational Mech. Anal., 29 (1968),
  pp.~32--57.

\bibitem{rayleigh}
{\sc L.~Rayleigh}, {\em On convection currents in a horizontal layer of fluid,
  when the higher temperature is on the under side}, Phil. Mag., 32 (1916),
  pp.~529--46.

\bibitem{salby}
{\sc M.~L. Salby}, {\em Fundamentals of Atmospheric Physics}, Academic Press,
  1996.

\bibitem{temam}
{\sc R.~Temam}, {\em Infinite-dimensional dynamical systems in mechanics and
  physics}, vol.~68 of Applied Mathematical Sciences, Springer-Verlag, New
  York, second~ed., 1997.

\bibitem{yudovich67a}
{\sc V.~I. Yudovich}, {\em Free convection and bifurcation}, J. Appl. Math.
  Mech., 31 (1967), pp.~103--114.

\bibitem{yudovich67b}
\leavevmode\vrule height 2pt depth -1.6pt width 23pt, {\em Stability of
  convection flows}, J. Appl. Math. Mech., 31 (1967), pp.~272--281.

\end{thebibliography}
\end{document}